\documentclass[amsthm]{elsart}
\usepackage{amsmath}
\usepackage{amsfonts}
\usepackage{amssymb}
\usepackage{amscd}
\usepackage{mathrsfs}
\usepackage{yjsco}
\usepackage{natbib}

\newcommand{\C}{\mathbb{C}}
\newcommand{\N}{\mathbb{N}}
\newcommand{\Z}{\mathbb{Z}}
\newcommand{\R}{\mathbb{R}}
\newcommand{\Q}{\mathbb{Q}}
\newtheorem{example}[thm]{Example}
\newtheorem{remark}[thm]{Remark}

\begin{document}
\begin{frontmatter}

\title{Toric Dynamical Systems}

\thanks{GC was supported by was supported by the
NSF (DMS-0553687) and DOE (BACTER Institute, DE-FG02-04ER25627).
AD was supported by UBACYT X042, CONICET PIP 5617 and ANPCyT PICT
20569, Argentina.
AS was supported by a Lucent Technologies Bell Labs Graduate Research Fellowship.
BS was supported by the NSF (DMS-0456960) and
DARPA ({\em Fundamental Laws in Biology}, HR0011-05-1-0057).}

\author{Gheorghe Craciun}
\address{Dept. of Mathematics, University of Wisconsin,
 Madison, WI 53706-1388, USA}
\ead{craciun@math.wisc.edu}

\author{Alicia Dickenstein}
\address{Dep.~de Matem\'atica,
FCEN, Universidad de Buenos Aires,
(1428), Argentina}
\ead{alidick@dm.uba.ar}

\author{Anne Shiu}
\address{Dept. of Mathematics, University of California, Berkeley, CA 94720-3840, USA}
\ead{annejls@math.berkeley.edu}

\author{Bernd Sturmfels}
\address{Dept. of Mathematics, University of California, Berkeley, CA 94720-3840, USA}
\ead{bernd@math.berkeley.edu}

%: ABSTRACT
\begin{abstract} Toric dynamical systems are known as complex balancing 
mass action systems in the mathematical chemistry literature, where 
many of their remarkable properties  have been established.
They include as special cases 
all deficiency zero systems and all detailed balancing systems.
One feature is that the steady state locus of a toric dynamical
system is a toric variety, which has a unique point within each
invariant polyhedron.
We develop the basic theory of toric dynamical systems in the context
of computational algebraic geometry and show that the associated moduli
space is also a toric variety.
It is conjectured that the complex balancing state is a global attractor.
We prove this for detailed balancing systems whose invariant 
polyhedron is two-dimensional and bounded.
 
\bigskip
\noindent {\em This paper is dedicated to the memory of Karin Gatermann (1961--2005).}
\end{abstract}


\end{frontmatter}

\section{Introduction} \label{sec:intro}

Toric dynamical systems describe mass-action kinetics with complex 
balancing states.  These systems have been studied extensively 
in mathematical chemistry, starting with the work of 
\citet{HornJackson72}, \citet{Fein72} and \citet{Horn72, Horn73}, 
and continuing with the deficiency theory in \citep{Fein79, FeinZero, 
FeinDet, Feinberg_1995}.
Mass-action kinetics has a wide range of applications in the physical 
sciences, and now it is beginning to play a role  in
systems biology \citep{CTF06, Gnad07, Guna, Sontag01}. Important special cases 
of these dynamical systems include recombination 
equations in population genetics  \citep{Akin79} and quadratic 
dynamical systems in  computer science \citep{RSW92}.

Karin Gatermann introduced the connection 
between mass-action kinetics and computational algebra.
Our work drew inspiration both from her publications \citep{Karin01, Karin02, Karin05}
and from her unpublished research notes on toric dynamical systems.
We wholeheartedly agree with her view that
{\em ``the advantages of toric varieties are well-known''}
\citep[page 5]{Karin01}.

We now review the basic set-up.
A {\em chemical reaction network} is a finite directed graph 
whose vertices are labeled by monomials and whose edges are labeled by parameters. 
The digraph is denoted $G = (V,E)$, with vertex set $V = \{1,2,\ldots,n\}$
and edge set $\,E \subseteq \{(i,j) \in V \times V : \,i\not= j \}$.
The node $i$ of $G$ represents the $i$th chemical complex and is labeled with the monomial
$$ c^{y_i} \,\,\, = \,\,\, c_1^{y_{i1}} c_2^{y_{i2}} \cdots  c_s^{y_{is}}. $$
Here $Y=(y_{ij})$ is an $n \times s$-matrix of non-negative integers.
The unknowns $c_1,c_2,\ldots,c_s$ represent the
concentrations of the $s$ species in the network,
and we regard them as functions $c_i(t)$ of time $t$.
The monomial labels are the entries in the row vector 
$$ \Psi(c) \quad = \quad \bigl( c^{y_1}, \, c^{y_2} , \, \ldots, \,c^{y_n} \bigr). $$
Each directed edge $(i,j) \in E$ is labeled by a positive parameter $\kappa_{ij}$ which
represents the rate constant in the reaction from 
the $i$-th chemical complex to the $j$-th chemical complex.
Note that if there is an edge from $i$ to $j$ and an edge from $j$
to $i$ then we have two unknowns $\kappa_{ij}$ and $\kappa_{ji}$.
Let $A_\kappa$ denote the negative of the {\em Laplacian} of
the digraph $G$. Hence $A_\kappa$ is the $n \times n$-matrix
whose off-diagonal entries
are the $\kappa_{ij}$ and whose row sums are zero.
Mass-action kinetics specified by the digraph $G$ is
the dynamical system
\begin{equation}
\label{CRN}
 \frac{d c}{dt} \quad = \quad \Psi(c) \cdot A_\kappa \cdot Y. 
 \end{equation}

A {\em toric dynamical system} is a dynamical system (\ref{CRN})
for which the algebraic equations
 $\,\Psi(c) \cdot A_\kappa  = 0\,$
admit a strictly positive solution
$c^* \in \R^s_{>0}$. Such a solution $c^*$ is a {\em steady
state} of the system, i.e., the $s$ coordinates of
$ \Psi(c^*) \cdot A_\kappa \cdot Y$ vanish.  The requirement
that all $n$ coordinates of $\,\Psi(c^*) \cdot A_\kappa  \,$ be zero
is stronger. The first to study toric dynamical systems, \citet{HornJackson72},
called these systems {\em complex balancing mass action systems} and
called $c^*$ a {\em complex balancing steady state}.
 A system (\ref{CRN}) being 
complex balancing (i.e., toric) depends on both the digraph G and the rate constants 
$\kappa_{ij}$.

\begin{example} \label{Ex:K3}
Let $s=2$, $n = 3$ and  let $G$ be the complete bidirected graph 
on three nodes labeled by
$c_1^2$, $c_1 c_2$ and $c_2^2$. Here the mass-action kinetics
system (\ref{CRN}) equals
\begin{equation}
\label{specialCRN}
\frac{d}{dt} \bigl(\, c_1,\,c_2\, \bigr) \quad = \quad
\begin{pmatrix} c_1^2 \,&\, c_1 c_2 \,&\, c_2^2 \end{pmatrix}
\cdot 
\begin{pmatrix}
 - \kappa_{12} - \kappa_{13} & \kappa_{12} & \kappa_{13} \\
  \kappa_{21} & -\kappa_{21} - \kappa_{23} & \kappa_{23} \\
  \kappa_{31} & \kappa_{32} & -\kappa_{31} - \kappa_{32}
  \end{pmatrix} \cdot
  \begin{pmatrix}
 \, 2 \,& \,0 \,\\
 \, 1 \,& \,1 \,\\
 \, 0 \,& \,2 \,\end{pmatrix}
  \end{equation}
  This is a toric dynamical system if and only if the following algebraic identity holds: 
  \begin{equation}
  \label{moduli1} 
  (\kappa_{21} \kappa_{31}+\kappa_{32} \kappa_{21}+\kappa_{23} \kappa_{31})
  (\kappa_{13} \kappa_{23}+\kappa_{21} \kappa_{13}+\kappa_{12} \kappa_{23})
\,   =\, (\kappa_{12} \kappa_{32}+\kappa_{13} \kappa_{32}+\kappa_{31} \kappa_{12})^2 .
   \end{equation}
 The equation (\ref{moduli1}) appears in \cite[Equation (3.12)]{Horn73}
 where it is derived from the necessary and sufficient 
 conditions for complex balancing in mass-action kinetics 
 given by \citet{Horn72}.
Our results in  Section~\ref{sec:ivch} provide a refinement of these conditions.
 
 Let us now replace $G$ by the digraph with four edges
 $(1,3), (2,1),(2,3),(3,1)$. This corresponds to setting
 $\kappa_{12} = \kappa_{32} = 0$ in (\ref{moduli1}).
  We can check that,  for this new $G$, the system  (\ref{CRN})
is not toric for any positive rate constants. Note that $G$ is not strongly connected. 
\qed
\end{example}

\smallskip

Among all chemical reaction networks, toric dynamical systems
have very remarkable properties. Some of these properties are
explained in \citep{Fein79}, starting with Proposition 5.3; 
see also \citep[Theorem 6.4]{Guna}. We shall review
them in detail in Sections~\ref{sec:ivch} and~\ref{sec:GAC}.  
From our point of view, the foremost among these
remarkable properties is that the
set $Z$ of all steady states is a toric variety 
\citep[\S 3]{Karin01}. Each trajectory of (\ref{CRN})
is confined to a certain {\em invariant polyhedron},
known to chemists as the {\em stoichiometric compatibility class},
which intersects the toric variety $Z$ in
precisely one point $c^*$.
In order to highlight the parallels between
toric dynamical systems and {\em toric models} in algebraic statistics
\citep[\S 1.2]{ASCB}, we shall refer to the steady state $c^*$ as the
{\em Birch point}; see \cite[Theorem 8.20]{Stu}.
In Example~\ref{Ex:K3}, the steady state variety $Z$ is a line through the origin, and
 the Birch point equals
$$ c^* \,\, = \,\, {\rm const} \cdot 
\bigl(\kappa_{12} \kappa_{32}+\kappa_{13} \kappa_{32}+\kappa_{31} \kappa_{12},\,
\kappa_{13} \kappa_{23}+\kappa_{21} \kappa_{13}+\kappa_{12} \kappa_{23} \bigr) .$$
Here the constant is determined because $c_1 + c_2$ is conserved along trajectories 
of (\ref{specialCRN}).

\smallskip

This article is organized as follows. In Section~\ref{sec:ivch} we develop 
the basic theory of toric dynamical systems within the context
of computational algebraic geometry.
 For each directed graph $G$ we introduce the moduli space
of toric dynamical systems on $G$.  This space
parametrizes  all rate constants $\kappa$
for which (\ref{CRN}) is toric. In Example \ref{Ex:K3}
this space is the hypersurface (\ref{moduli1}).  Theorem \ref{cayley}
states that this moduli space is itself a toric variety in a suitable
system of coordinates. These coordinates are the maximal non-zero minors of
the Laplacian of  $G$, and their explicit form as positive polynomials in the $\kappa_{ij}$
is given by the {\em Matrix-Tree Theorem} \citep[\S 5.6]{Sta}.
Our results in Section~\ref{sec:ivch} furnish a two-fold justification
for attaching the adjective ``toric'' to chemical reaction
 networks with complex balancing, namely,
 both the steady state variety and the moduli space are toric.
In addition, the subvariety of reaction networks with detailed balancing is toric.

In Section~\ref{sec:GAC} we introduce the {\em Global Attractor Conjecture}
which states that the Birch point is a global attractor for any toric dynamical system.
More precisely, we conjecture that
all trajectories beginning at strictly positive vectors $c^0$ will converge to the
Birch point $c^*$ in the invariant polyhedron of $c^0$. The conjecture
is currently open, even for {\em deficiency zero systems} (cf. Theorem~\ref{cayley}).
\citet{DAS07} found a proof for
a class of  ``monotone'' deficiency zero networks where the
monomials $c^{y_i}$ involve distinct unknowns.
We prove the conjecture in Section~\ref{sec:GACresults} 
for toric dynamical systems with detailed balancing that evolve in a bounded 
polygon in s-dimensional space. The algebraic theory of
detailed balancing systems is developed in Section~\ref{sec:detailed}.

\section{Ideals, Varieties and Chemistry} \label{sec:ivch}

This section concerns the connection between chemical reaction network theory
and toric geometry. We use the language of ideals and varieties
as in \citep{CLO}. Our reference on toric geometry and its
relations with computational algebra is \citep{Stu}.
With regard to the dynamical system (\ref{CRN}),
we use the notation from \citep[\S 5]{Fein79} and \cite[\S 3]{Guna}
which has the virtue of separating the roles
played by the concentrations $c_i$, the monomials $c^{y_{i}}$, and
the rate constants $\kappa_{ij}$.

To study the dynamical system (\ref{CRN}) algebraically, we work in the polynomial ring
$$ \Q[c,\kappa] \,\, = \,\, \Q \bigl [ \{c_1,c_2,\ldots,c_s \} \,\cup \,\{ \kappa_{ij}: (i,j) \in E\}  \bigr], $$
and we introduce various ideals in this polynomial ring.
First, there is the {\em steady state ideal}
$\,\langle \Psi(c) \cdot A_\kappa \cdot Y \rangle\,$
which is generated by the $s$ entries of the
row vector on the right hand side of (\ref{CRN}).
Second, we consider the ideal
$\,\langle \Psi(c) \cdot A_\kappa \rangle \,$ which is 
generated by the $n$ entries of the row vector
$\,\Psi(c) \cdot A_\kappa$. The generators
of both ideals  are linear in the $\kappa_{ij}$
but they are usually non-linear in the $c_i$.
Next, we define the {\em complex balancing ideal} of $G$ to be the
following ideal quotient whose
generators are usually non-linear in the $\kappa_{ij}$:
$$ C_G \quad := \quad   
\bigl( \,\langle \Psi(c) \cdot A_\kappa  \rangle \,: \,
(c_1 c_2 \cdots c_s)^\infty \bigr). $$
We have thus introduced three ideals in $\Q[c,\kappa].$
They are related by the inclusions
$$ \langle \Psi(c) \cdot A_\kappa \cdot Y \rangle
\,\, \subseteq \,\, \langle \Psi(c) \cdot A_\kappa \rangle \,\,\,
\subseteq \,\,\, C_G. $$

If $I$ is any polynomial ideal then we write $V(I)$ for its complex variety.
Likewise, we define the positive variety $V_{>0}(I)$
and the non-negative variety $V_{\geq 0}(I)$. They consist 
of all points in $V(I)$ whose coordinates are real and
positive or, respectively, non-negative. Our algebraic approach to
chemical reaction network theory focuses on the study of 
these varieties. The inclusions of ideals above imply the following inclusions of varieties:
\begin{equation}
\label{inclusions}
V(C_G) \,\,\,\,\subseteq\,\,\,\, V \bigl( \langle \Psi(c) \cdot A_\kappa \rangle \bigr) \, \,\, \, \subseteq \,\,\,\,
V \bigl( \langle \Psi(c) \cdot A_\kappa \cdot Y \rangle \bigr). \,\,\,
\end{equation}
The definition of $C_G$ by means of saturation implies that
the left hand inclusion becomes equality when 
we restrict to the points with all coordinates non-zero.
In particular,
\begin{equation}
\label{allequal}
 \,\, V_{> 0}(C_G) \,\,\, = \,\,\, V_{> 0} \bigl( \langle \Psi(c) \cdot A_\kappa \rangle \bigr). \qquad
\end{equation}

Recall from \citep{Stu} that a {\em toric ideal} is a prime ideal which is generated by binomials.
We soon  will replace $C_G$  by a subideal $T_G$ which is toric.
This is possible by Proposition 5.3 (ii,iv) in \citep{Fein79} or Theorem 6.4 (3) in \citep{Guna}, which 
essentially state that $V_{>0}(C_G)$ is a positive toric variety.
But let us first examine the case when 
$C_G$ is a toric ideal already.

\begin{example} \label{ex:lawrence}
Suppose that each chemical complex appears in only 
one reaction, and each reaction is bi-directional.
Hence $n =2m$ is even and, after relabeling, we have $\,
E = \{(1,2),(2,1), (3,4), (4,3), \ldots,(n{-}1,n),(n,n{-}1)\}$. 
We start with the binomial ideal
$$
\langle \Psi(c) \cdot A_\kappa \rangle \, \, = \,\,
\bigl\langle \,
\kappa_{12} c^{y_1} - 
\kappa_{21} c^{y_2},\,
\kappa_{34} c^{y_3} - 
\kappa_{43} c^{y_4}, \,\ldots \,,\,
\kappa_{n-1,n} c^{y_{n-1}} - \kappa_{n,n-1} c^{y_n} \bigr\rangle .
$$
The complex balancing ideal $C_G$ is a saturation of 
$\,\langle \Psi(c) \cdot A_\kappa \rangle $, and it coincides with the 
toric ideal of the extended Cayley matrix in the proof of Theorem \ref{moduli}.
There are many programs for computing toric ideals.
For instance, the methods in \cite[\S 12.A]{Stu} are
available in {\tt maple} under the command 
{\tt ToricIdealBasis}. Explicitly, the complex balancing ideal
 $C_G$ is generated by all binomials
$\, \kappa^{u_+} c^{v_+} \,-\, \kappa^{u_-} c^{v_-} \,$
where 
\begin{equation}
\label{linsys}
\sum_{i=1}^m
u_{2i-1,2i} (y_{2i-1}-y_{2i}) \,=\, v 
\quad {\rm and} \quad
u_{2i-1,2i} + u_{2i,2i-1} = 0 \,\,\hbox{for $i=1,2,\ldots,m$}.
\end{equation}
Eliminating  $c_1,\ldots,c_s$
from $C_G$, we obtain the ideal of all binomials 
$\,\kappa^{u_+} - \kappa^{u_-}\,$
where $u \in \N^E$ satisfies (\ref{linsys}) with $v  = 0$.
This is the moduli ideal $M_G$ to be featured in
Theorems \ref{moduli} and \ref{cayley} below. It is 
a prime binomial ideal of Lawrence type \citep[\S 7]{Stu}.
\qed
\end{example}

\smallskip

Let us next assume that $G = (V,E)$ is an arbitrary digraph
with $n$ nodes which is {\em strongly connected}. This means
that, for any two nodes $i$ and $j$, there exists a directed
path from $i$ to $j$. In this case the
matrix $A_\kappa$ has rank $n-1$,
and all its minors of size  $(n-1) \times (n-1)$ are non-zero.
The next result gives a formula for these comaximal minors.

Consider any directed subgraph $T$ of $G$
whose underlying graph is a tree.
This means that $T$ has $n-1$ edges
and contains no cycle. We write
$\kappa^T$ for the product of all
edge labels of the edges in $T$. This is
a squarefree monomial in $\Q[\kappa]$.
Let $i$ be one of the nodes of $G$.
The directed tree $T$ is called an {\em $i$-tree} if 
the node $i$ is its unique sink, i.e., all edges are directed
towards node $i$. We introduce the following polynomial
of degree $n-1$:
\begin{equation}
\label{Kpolynomial}
 K_i \,\,\, = \,\, \sum_{T \,\text{an $i$-tree}} \kappa^T .
 \end{equation}
The following result is a restatement of
the {\em Matrix-Tree Theorem} \citep[\S 5.6]{Sta}.
 
\begin{prop}
\label{MTT}
Consider a submatrix of $A_\kappa$
obtained by deleting the $i^{th}$ row and any one of the columns.
The signed determinant of this 
$(n{-}1) \times (n{-}1)$-matrix equals $(-1)^{n-1} K_i$.
\end{prop}

\smallskip

This minor is independent of the choice of columns
because the row sums of $A_\kappa$ are zero.
Combining Proposition \ref{MTT} with a little
linear algebra leads to the following corollary:

\begin{cor}
The complex balancing ideal
$C_G$ contains the polynomials
$\, K_i c^{y_j} - K_j c^{y_i}$.
\end{cor}

\smallskip

We now form the ideal
generated by these $\binom{n}{2}$ polynomials
and we again saturate with respect to $c_1 c_2 \cdots c_s$.
The resulting ideal $T_G$ will be called the {\em toric balancing ideal}:
$$ T_G \quad := \quad \bigl( \, \langle
 K_i c^{y_j} - K_j c^{y_i} \,:\,
 1 \leq i < j \leq n \rangle \,: \,
(c_1 c_2 \cdots c_s)^\infty \bigr). $$
It is thus  natural to consider $T_G$ as an ideal in the
polynomial subring
$$ \Q[c,K] \,\, = \,\,\Q[c_1,\ldots,c_s,K_1,\ldots,K_n]
\quad \subset\,\,\,\, \Q[c,\kappa]. $$
The claim that this is a polynomial ring is the content of the
following lemma.

\begin{lem}
The polynomials $K_1,\ldots,K_n \in \Q[\kappa]$ are
algebraically independent over~$\Q$.
\end{lem}
\begin{pf}
Let $K_i' \in \Q [\kappa_1, \kappa_2, \dots, \kappa_n ]$ denote
the polynomial obtained from $K_i$ by substituting the new unknown $\kappa_i$ for all $\kappa_{ij}$.  We need only verify that the $K'_i$ are algebraically independent, because an algebraic relation among the $K_i$ would be satisfied by the $K_i'$ as well.  Our polynomials are 
$$ K_i' \quad = \quad \text{(number of $i$-trees in $G$)} \cdot \prod_{t \neq i} \kappa_t.$$  
The $n$ squarefree monomials 
$ \prod_{t \neq i} \kappa_t$ (for $i=1\dots n$)
are algebraically independent
because an algebraic dependence among these monomials would specify a dependence among $1/{\kappa_1}, 1/{\kappa_2}, \dots, 1/{\kappa_n}$.  Hence, $K_1', K_2', \dots, K_n'$ are algebraically independent.
\end{pf}
\smallskip

We new discuss the toric balancing ideal $T_G$.

\begin{prop} 
\label{prop:itstoric}
The toric balancing ideal $T_G$ is a toric ideal in
$\Q[c,K]$. Moreover, the ideal $T_G$ is generated by the binomials
$\,K^{u_+} \cdot c^{(uY)_-}  - K^{u_-} \cdot c^{(uY)_+}\,$
where $u $ is any row vector in $\Z^n$
whose coordinate sum $u_1+u_2 + \cdots +u_n $ is equal to zero.
\end{prop}

\begin{pf}
Let $\Delta$ denote the edge-node incidence matrix of the complete directed
graph on $n$ nodes. Thus $\Delta$ is the $\binom{n}{2} \times n$-matrix
whose rows are $e_i-e_j$ for $1 \leq i < j \leq n$. We also consider the
$n \times (n+s) $-matrix $\,\bigl(\,-Y \,\,\, {\bf I}_n \,\bigr) $.
The binomials $K_i c^{y_j} - K_j c^{y_i}$ which define the ideal $T_G$
correspond to the rows of the $\binom{n}{2} \times (n+s)$-matrix
$\, \Delta \cdot \bigl(\,-Y \,\,\, {\bf I}_n \,\bigr) $, and
the binomial $\,K^{u_+} \cdot c^{(uY)_-}  - K^{u_-} \cdot c^{(uY)_+}\,$
corresponds to the row vector $\, U \cdot \Delta \cdot \bigl(\,-Y \,\,\, {\bf I}_n \,\bigr) $,
where $U$ is any row vector of length $\binom{n}{2}$ such that $u = U \cdot \Delta$.
The binomial $\,K^{u_+} \cdot c^{(uY)_-}  - K^{u_-} \cdot c^{(uY)_+}\,$
is a $\Q[c_1^{\pm 1}, \ldots, c_s^{\pm 1}, K_1,\ldots,K_n]$-linear
combination of the binomials $K_i c^{y_j} - K_j c^{y_i}$.
This shows that $T_G$ is the {\em lattice ideal}
in $\Q[c,K]$ associated with the lattice spanned by the rows of
$\, \Delta \cdot \bigl(\,-Y \,\,\, {\bf I}_n \,\bigr) $, i.e., there
are no monomial zero-divisors
modulo $T_G$. To see that $T_G$ is actually a toric ideal, i.e.~$T_G$
is prime, it suffices to note that $\Z^{n+s}$ modulo the lattice spanned
by the rows of $\, \Delta \cdot \bigl(\,-Y \,\,\, {\bf I}_n \,\bigr) \,$ is free abelian
of rank $s+1$. Indeed, the latter matrix has rank $n-1$, and its
$(n-1)\times (n-1)$-minors span the unit ideal in the ring of integers  $\Z$,
because each $(n-1) \times (n-1)$-minor of $\Delta$ is either $+1$ or $-1$. 
\end{pf}

The variety of $T_G$ is a toric variety in ${\rm Spec}\,\Q[c,K]$, 
but we continue to regard it as a subvariety of $\C^{s} \times \C^E\,$
(or of ${\rm Spec}\,\Q[c,\kappa]$). In this interpretation we have
\begin{equation}
\label{allequal2}
V_{>0}(T_G) \,\,\, = \,\,\,
 \,\, V_{> 0}(C_G) \,\,\, = \,\,\, V_{> 0} \bigl( \langle \Psi(c) \cdot A_\kappa \rangle \bigr)   .
\end{equation}
Thus $T_G$ still correctly describes the steady state locus
of the toric dynamical system.
The equation (\ref{allequal2}) holds because
the matrix $A_\kappa$ has rank $n-1$
over the rational function field $\Q(\kappa)$, and
the vector $\,(K_1,K_2,\ldots,K_n)\,$ spans
its kernel under left multiplication.

Finally, the following elimination ideal
is called the {\em moduli ideal} of the digraph $G$:
\begin{equation}
\label{defMG}
 M_G \quad = \quad    T_G \,\,\cap \,\, \Q[\kappa]  .
 \end{equation}
Here $\Q[\kappa]$ is the polynomial ring in only
the edge unknowns $\kappa_{ij}$. The generators of $M_G$
are obtained from the generators of $C_G$ by eliminating
the unknown concentrations $c_i$.
For instance, if $G$ is the complete bidirected
graph on $c_1^2$, $c_1 c_2$ and $c_2^2$
as in Example \ref{Ex:K3} then the moduli ideal
$M_G$ is the principal ideal generated
by $K_1 K_3 - K_2^2$. This coincides
with condition (\ref{moduli1}) because
$\,K_1 = \kappa_{21} \kappa_{31} + \kappa_{32} \kappa_{21} + \kappa_{23} \kappa_{31}$,
and similarly for $K_2$, $ K_3$.

Suppose now that $G$ is an arbitrary directed graph,
and let $l$ be the number of connected components of $G$.
If one of the components $G_i$ fails to be strongly connected,
then $V_{>0}(C_{G_i})$ is empty  and hence $V_{>0}(C_G)$ is empty, by \citep[Remark 5.2]{Fein79}.
In that case we define $T_G$ and $M_G$ to
be the ideal generated by $1$. If each connected component $G_i$ of $G$
is strongly connected then we define the toric steady state ideal as
$$ T_G \quad := \quad \bigl( \,
(\,T_{G_1} + T_{G_2} + \cdots + T_{G_l}): (c_1 c_2 \cdots c_s)^\infty \bigr). $$
The moduli ideal $M_G$ is defined as before in (\ref{defMG}).
The equality in (\ref{allequal2}) still
holds and this positive variety is in fact non-empty.
Here is the first main result of this section:

\begin{thm} \label{moduli}
The equations (\ref{CRN}) specify
a toric dynamical system if and only if
the positive vector of
rate constants  $\kappa_{ij}$
lies in the toric variety $V(M_G)$.
In this case, the set of steady states of (\ref{CRN}) with all $c_i > 0$
equals the set of positive points on the toric variety $V(T_G)$.
\end{thm}

\begin{pf}
The positive variety $V_{>0}(T_G)$ consists of all pairs $(c,\kappa)$ where $\kappa$ is 
a strictly positive vector of rate constants and $c$ is a strictly positive solution of the complex 
balancing equations $\Psi(c) \cdot A_\kappa = 0$. The elimination in (\ref{defMG}) corresponds 
to the map of toric varieties $\, V(T_G) \,\rightarrow \, V(M_G)\,$ given by 
$\,(c,\kappa) \,\mapsto \, \kappa$. This map is a dominant morphism (by definition of $M_G$), 
so its image is Zariski dense in $V(M_G)$. The restriction to real positive points,
$\, V_{>0}(T_G) \rightarrow  V_{>0}(M_G)$,  is a homomorphism of abelian groups $\,(\R_{>0})^*\,$ 
whose image is dense, so it is the monomial map specified by a matrix with maximal row rank.
It follows that this restriction is  surjective, and this proves our first assertion.
The second assertion follows from \citep[Proposition 5.3]{Fein79}.
\end{pf}

We now justify calling $V(M_G)$ a toric variety
by writing $M_G$ explicitly
as a toric ideal in $\Q[K]$. As before, $G$ is a directed graph with $n$ nodes
labeled by monomials $c^{y_1}, \ldots , c^{y_n}$. We assume that
each connected component $G_1,G_2,\ldots,G_l$ 
of $G$ is strongly connected,
for otherwise $M_G \,=\, \langle 1 \rangle$.
Let $Y_i$ denote the matrix with $s$ rows whose
columns are the vectors $y_j$ where $j$
runs over the nodes of the component $G_i$.
We define the {\em Cayley matrix} 
$$ {\rm Cay}_G(Y) \quad = \quad
\begin{pmatrix}
  \,\, Y_1 \,\,& \,\,Y_2 \,\,&  \,\,\cdots \,\,&  \,\,Y_l \,\, \\
{\bf 1} & {\bf 0} & \cdots & {\bf 0} \\
{\bf 0} & {\bf 1} &  \cdots & {\bf 0} \\
\vdots & \vdots & \ddots & \vdots \\
{\bf 0} & {\bf 0} &  \cdots & {\bf 1}
\end{pmatrix}.
$$
This is an $(s + l) \times n$-matrix.
Here ${\bf 1}$ and ${\bf 0}$
are appropriate row vectors with all
entries $1$ and $0$ respectively.
The term ``Cayley matrix'' comes from
geometric combinatorics, and it refers to the
Cayley trick in elimination theory 
\citep{HRS00}.

Let $S $ denote the linear subspace of $\R^s$ which is spanned
by the {\em reaction vectors} $y_j-y_i$ where $(i,j) \in E$.
This space is known in chemistry as the
{\em stoichiometric subspace}.
We write $\sigma = {\rm dim}(S)$ for its dimension.
The quantity $\,\delta \,:=\,n-\sigma -l$ is known as the
{\em deficiency} of the chemical reaction network $G$. 
 For instance,  $\delta = 3-1-1 = 1$  in Example \ref{Ex:K3}.

\begin{remark} \label{rankofcayley}
The rank of the Cayley matrix ${\rm Cay}_G(Y)$ equals
$\sigma + l$. Hence the deficiency $\delta$
of the reaction network coincides with the
dimension of the kernel of the Cayley matrix.
\end{remark}

\smallskip

The following theorem is the second main result in this section.

\begin{thm} \label{cayley}
The moduli ideal $M_G$ equals the toric ideal 
of the Cayley matrix ${\rm Cay}_G(Y)$, i.e.~$M_G$
is the ideal in $\Q[K]$ generated by
all binomials $\,K^u - K^v\,$ where $u,v \in \N^n$
satisfy ${\rm Cay}_G(Y) \cdot (u-v) = 0 $.
The codimension of this toric ideal equals the deficiency $\delta$.
\end{thm}

\begin{pf}
Let ${\bf Id}_s$ denote the $s \times s$ identity matrix
and consider the extended Cayley matrix
$$
\begin{pmatrix}
\,\, - {\bf Id}_s \,\,\, & \,\, Y_1 \,\,& \,\,Y_2 \,\,&  \,\,\cdots \,\,&  \,\,Y_l \,\, \\
{\bf 0} &   {\bf 1} & {\bf 0} & \cdots & {\bf 0} \\
{\bf 0} &  {\bf 0} & {\bf 1} &  \cdots & {\bf 0} \\
{\bf 0} &  \vdots & \vdots & \ddots & \vdots \\
{\bf 0} &  {\bf 0} & {\bf 0} &  \cdots & {\bf 1}
\end{pmatrix}.
$$
The toric ideal of this matrix is precisely the
toric balancing ideal $T_G$, where  the unknowns
$c_1,c_2,\ldots,c_s$ correspond to the first $s$
columns. Deleting these $s$ columns corresponds to
forming the elimination ideal $M_G$ as in (\ref{defMG}).
This shows that $M_G$ is the toric ideal of the matrix
${\rm Cay}_G(Y)$.
The dimension of the affine toric variety $\,V(M_G) \,$ in $\, \C^n\,$
is equal to $\,\sigma+l = {\rm rank}({\rm Cay}_G(Y))$,
and hence its codimension equals the deficiency
$\,\delta = n - \sigma - l$.  
\end{pf}

We conclude that $\,V_{>0}(M_G)\,$ is a 
positive toric variety of codimension $\delta$ in $\R_{>0}^n$.
The moment map of toric geometry establishes a 
natural bijection between $\,V_{>0}(M_G)$ and the
interior of the {\em Cayley polytope}, which is
the convex hull of the columns of ${\rm Cay}_G(Y)$.

In summary, given any chemical reaction network whose components
are strongly connected, we have shown that the positive
toric variety of the Cayley polytope equals the moduli space
$V_{>0}(M_G)$ of toric dynamical systems on $G$.
The deficiency $\delta$  is precisely the codimension of this moduli space.
In particular, if the deficiency is zero
 then the Cayley polytope is a simplex
and (\ref{CRN}) is toric for all rate constants $\kappa_{ij}$.
Moreover, the positive steady states of a toric dynamical system
form a positive toric variety $V_{>0}(T_G)$.

\section{The Global Attractor Conjecture and Some Biological Applications} \label{sec:GAC}

We now consider a fixed toric dynamical system
or, equivalently, a chemical reaction network (\ref{CRN})
 which admits a complex balancing state.
The underlying directed graph
 $G = (V,E)$ has $n$ nodes labeled by monomials $c^{y_1}$, $c^{y_2}$, ..., $c^{y_n}$,
 and we specify positive rate constants by fixing a  
 point $\kappa^0$ in the moduli space $V_{>0}(M_G)$.
We also fix a strictly positive vector $c^0 \in \R_{>0}^s$
which represents the initial concentrations of the
$s$ species. The equations (\ref{CRN}) describe
the evolution of these concentrations over time. We seek
to understand the long-term behavior of the trajectory
which starts at $c^0$, that is, $c(0) = c^0$.

Let $\,T_G(\kappa^0)\,$ denote the toric ideal in $\R[c]$
obtained from $T_G$ by substituting the specific rate constants
$\,\kappa^0_{ij} \in \R_{>0}\,$ for the unknowns $\kappa_{ij}$.
Then $V_{>0}(T_G(\kappa^0))$ coincides with the set of all
steady states of the toric dynamical system (\ref{CRN}).
The following result is well-known:

\begin{prop}
\label{Birch} {\rm [Existence and Uniqueness of the Birch Point] }
The affine subspace $\,c^0+S\,$ of $\R^s$
intersects the  positive toric variety $V_{>0}(T_G(\kappa^0))$
in precisely one point $c^*$.
\end{prop}

\smallskip

For a proof and references in the chemistry literature see \citet{HornJackson72}; 
a different proof can be found in \citet[Proposition 5.3]{Fein79} or
\citet[Proposition 6.4]{Guna}. We remark that variants of Proposition
\ref{Birch} are ubiquitous across the mathematical
sciences, and the result has been rediscovered many times.
In statistics, this result is known as Birch's Theorem; see
\cite[Theorem 1.10]{ASCB}. To stress the link with
toric models in algebraic statistics we call $c^*$ the {\em Birch point} of
the toric dynamical system (\ref{CRN}) with starting point~$c^0$.

The right hand side of (\ref{CRN}) is always a vector in the
stoichiometric subspace  $\,S = \R \{  y_j-y_i \,: (i,j) \in E \}$.
Hence the trajectory starting at $c^0$ stays in the affine subspace
$\,c^0 + S$. In fact, concentrations remain non-negative, so the
trajectory stays in  $\, P \, := \, \,(c^0 +S) \cap \R_{\geq 0}^s$.
We call $P$ the {\em invariant polyhedron}.
Chemists use the term {\em stoichiometric compatibility class}
for $P$. The relative interior of $P$ in  $c^0  + S$
is denoted by $\,P^o := (c^0 +S) \cap \R_{> 0}^s$.

\begin{prop}\label{prop:Birch}
The Birch point $\,c^*\,$ is the unique point in the invariant 
polyhedron $P$ for which the {\em transformed entropy function}
\begin{equation} \label{eq:entropy}
E(c) \,\,\, = \,\,\, \sum_{i=1}^s \bigl(\, c_i \cdot {\rm log}(c_i) \,
-\, c_i \cdot {\rm log}(c_i^*)  \,-\, c_i \, + c_i^*\bigr)
\end{equation}
is a {\em strict Lyapunov function} of the toric dynamical system (\ref{CRN}).
This means the following:
\begin{itemize}
\item[(a)] For all $c \in P$ we have $E(c) \geq 0$  and equality holds if and
only if $c=c^*$,
\item[(b)]  we have $\,d E (c)/dt \leq 0\,$ along any trajectory $c(t)$ in $P$, and
\item[(c)]  equality in (b)
holds at a point $t$ of any trajectory $c(t)$ in
$P^\circ$ if and only if $c(t) = c^*$.
\end{itemize}
\end{prop}

\smallskip

This proposition was proved by \citet{HornJackson72}. A different proof can be found in \citep{Fein79}; see especially Proposition 5.3 and its corollaries; see also \citep[Theorem 6.4]{Guna} 
and the paragraph before it. We suggest comparing this with the proof of \citep[Theorem 1.10]{ASCB}.

Any trajectory of the toric dynamical system (\ref{CRN}) which 
starts in the relatively open polyhedron $\,P^o = (c^0 +S) \cap \R_{> 0}^s\,$
will stay in the closed polyhedron $\,P =  (c^0 +S) \cap \R_{\geq 0}^s $; actually, 
it is not hard to show that $P^o$ is an invariant set.
The main conjecture below states that any such trajectory converges to the 
Birch point. This conjecture was first formulated by \citet{Horn74}. 
A steady state $x$ in $P^o$ is called a {\em global attractor} if any 
trajectory that begins in $P^o$ converges to~$x$.  

\medskip

\noindent{\bf Global Attractor Conjecture.} \emph{For any toric dynamical system 
(\ref{CRN}) and any starting point $c^0$, the Birch point $c^*$ is a global attractor of the invariant set 
$\,P^o = (c^0 +S) \cap \R_{>0}^s$.}

\medskip

An important subclass of toric dynamical systems consists of the
chemical reaction networks of deficiency zero.
If the deficiency $\delta = n - \sigma - l$ is zero then
the moduli ideal $M_G$ is the zero ideal,
by Theorem \ref{cayley}, and 
(\ref{CRN}) is toric for all choices
of rate constants. As remarked in the Introduction, the Global 
Attractor Conjecture is open even for deficiency zero systems.
Our last section is devoted to partial results on the
conjecture. First, however, we discuss  biological examples
which illustrate the concepts developed so far.

\begin{example} {\bf [Networks with trivial moduli]}
We expect that our toric approach will be useful for parametric
analyses of chemical reaction networks in {\em systems biology}.
Analyses of this kind include \citep{KSP07}, \citep{Gnad07} 
and \citep{Sontag01}.
Many of the explicit examples we found in the literature
have trivial toric moduli in the sense that either
$M_G$ is the unit ideal or $M_G$ is the zero ideal.

If  $M_G = \langle 1 \rangle $ then (\ref{CRN}) is {\bf never} a
toric dynamical system regardless of what the $\kappa_{ij}$ are. This happens 
when components of $G$ are not strongly connected. Examples include 
 {\em Michaelis-Menten kinetics} and the {\em covalent modification cycle} in 
 \cite[\S 5]{Guna}.
 If $M_G = \{0\}$ then the network has deficiency zero
 and (\ref{CRN}) is {\bf always} a toric dynamical system, regardless
 of what the $\kappa_{ij}$ are.
 Examples include the cycle in \citep[Equation (9)]{KSP07},
 the monotone networks in \citep{DAS07},
 and the following network which is taken from \citep{Gnad07}.
 
 The {\em ligand-receptor-antagonist-trap network}
 has $s=8$ species and $n=8$ complexes.
 This network $G$ has four reversible reactions
 which we write in binomial notation:
 \begin{equation}
 \label{trap}
 \kappa_{15} \cdot c_5 c_6 \,-\,\kappa_{51} \cdot c_1, \,\,
 \kappa_{26} \cdot c_6 c_7 \,-\,\kappa_{62} \cdot c_2,\,\,
 \kappa_{37} \cdot c_7 c_8 \,-\, \kappa_{73} \cdot c_3,\,\,
 \kappa_{48} \cdot c_8 c_5 \,-\, \kappa_{84} \cdot c_4.
 \end{equation}
Here $l = 4$ and $\sigma = 4$, so that $\delta = 0$.
In the algebraic notation of Section~\ref{sec:ivch}, the toric ideal $T_G$ 
equals the complex balancing ideal
$\,\langle \Psi(c) \cdot A_\kappa \rangle  \,$ and is generated by
the four binomials in (\ref{trap}).
 Eliminating $c_1,c_2,\ldots,c_8$ as prescribed by
(\ref{defMG}) yields the zero ideal $M_G = \{0\}$. \qed
 \end{example}

\begin{example} {\bf [DHFR catalysis]}
Here are some examples from systems biology which show a more complicated dynamical behaviour. 
We consider the reaction network in (Craciun, Tang and Feinberg, 2006, Figure 5);  this reaction network has several positive equilibria for some values of the  reaction rate parameters (see Craciun, Tang and Feinberg, 2006, Figure 7). %Therefore, it cannot be a toric dynamical system for all parameter values. 
This reaction network allows for inflow and outflow of some chemical species; in the language of deficiency theory, we say that one of the complexes of this reaction network is the zero complex (see \cite{Fein79}), i.e., one of the vectors $y_i$ is zero. 
Note that the group A of reactions in this network has almost the same structure as mechanism 6 in (Craciun, Tang and Feinberg, 2006, Table 1), shown below:
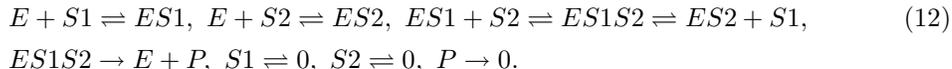
\begin{eqnarray}
\label{rn_two_substrate}
& & E+S1 \rightleftharpoons ES1,\ E+S2\rightleftharpoons ES2,\ ES1+S2 \rightleftharpoons ES1S2\rightleftharpoons ES2+S1,\ \\ \nonumber
& & ES1S2\to E+P,\ S1 \rightleftharpoons 0,\ S2 \rightleftharpoons 0,\ P \to 0. 
\end{eqnarray}
Like the more complicated DHFR catalysis network, the network (\ref{rn_two_substrate}) also has several positive equilibria for some values of the reaction rate parameters. It is easy to compute the deficiency of this simpler mechanism: the number of complexes is $n=12$ (including the zero complex), the number of linkage classes is $l=4$ (including the linkage class that contains the inflow and outflow reactions for the substrates $S1$, $S2$ and the product $P$), and the dimension of its stoichiometric subspace is $\sigma=6$. Therefore the deficiency of the network (\ref{rn_two_substrate}) is $\delta = 12-4-6 = 2$. This network cannot be toric for any choice of the
constant rates because it is not weakly reversible.
If we make all reactions reversible in (\ref{rn_two_substrate}), then the complexes, the linkage classes, and the stoichiometric subspace do not change, so the deficiency of the reversible version of (\ref{rn_two_substrate}) is also 2.
\end{example}

\begin{example} {\bf [Recombination on the $3$-cube]} \
In {\em population genetics} \citep{Akin79,Akin82},
the evolution of a population is  modeled
by a dynamical system whose right hand side is the sum of three terms, 
corresponding to {\em mutation}, {\em selection} and {\em recombination}.
The contribution made by recombination alone is a quadratic
dynamical system \citep{RSW92} which can be written in
the form (\ref{CRN}). In our view, toric 
dynamical systems are particularly well-suited
to model recombination. Here we consider a 
population of three-locus diploids, so the
underlying {\em genotope} of the haploid gametes is the
standard $3$-dimensional cube \citep[Example 3.9]{BPS07}.
The eight vertices of the cube are the genotypes.
They now play the role of the species in chemistry:
$$ \!\!\!\!\!\!\!\!\!\!\! s=8 \qquad \begin{matrix}
\,{\rm genotypes} &\quad & [000] & [001] & [010] & [011] & [100] & [101] & [110] & [111] \\
{\rm frequencies} &\quad & c_1 & c_2 & c_3 & c_4 & c_5 & c_6  & c_7  & c_8 .
\end{matrix}
$$
The {\em recombination network} $G$ has $n=16$ nodes
corresponding to the pairs of genotypes which are not
adjacent on the cube. There are twelve bidirectional edges, 
representing interactions, and we label them using the notation of
\citep[Example 3.9]{BPS07}.
Six of the interactions  correspond to {\em conditional epistasis}:
$$ \begin{matrix}
 [000] + [110] \leftrightarrow [010] + [100] 
 & \quad   \kappa_{1,2} \cdot c_1 c_7 - \kappa_{2,1} \cdot c_3 c_5  \quad
 & K_1 = \kappa_{2,1} \,\,{\rm and} \,\, K_2 = \kappa_{1,2} \\
 [001] + [111] \leftrightarrow [011] + [101] 
 & \quad   \kappa_{3,4} \cdot c_2 c_8 - \kappa_{4,3} \cdot c_4 c_6  \quad
 & K_3 = \kappa_{4,3} \,\,{\rm and} \,\, K_4 = \kappa_{3,4} \\
  [000] + [101] \leftrightarrow [001] + [100] 
 & \quad   \kappa_{5,6} \cdot c_1 c_6 - \kappa_{6,5} \cdot c_2 c_5  \quad
 & K_5 = \kappa_{6,5} \,\,{\rm and} \,\, K_6 = \kappa_{5,6} \\
  [010] + [111] \leftrightarrow [011] + [110] 
 & \quad   \kappa_{7,8} \cdot c_3 c_8 - \kappa_{8,7} \cdot c_4 c_7  \quad
 & K_7 = \kappa_{8,7} \,\,{\rm and} \,\, K_8 = \kappa_{7,8} \\
  [000] + [011] \leftrightarrow [001] + [010] 
 & \quad   \kappa_{9,10} \cdot c_1 c_4 - \kappa_{10,9} \cdot c_2 c_3  \quad
 & K_9 = \kappa_{10,9} \,\,{\rm and} \,\, K_{10} = \kappa_{9,10} \\
  [100] + [111] \leftrightarrow [101] + [110] 
 & \,\,\,   \kappa_{11,12} \cdot c_5 c_8 - \kappa_{12,11} \cdot c_6 c_7  \,\,\,
 & K_{11} = \kappa_{12,11} \,{\rm and} \, K_{12} = \kappa_{11,12} .
  \end{matrix}
  $$
  Secondly, we have {\em marginal epistasis}, giving rise to
the six pairwise interactions among
$$ \begin{matrix}
{\rm four} \, {\rm complexes} & \quad & [000]+[111] &\quad& [001]+[110] &\quad& [010]+[101] &\quad& [100]+[011] \\ %corrected typo in this line
\hbox{four monomials} & \quad & K_{13}\cdot c_1 c_8 && 
K_{14} \cdot c_2 c_7 && K_{15} \cdot c_3 c_6 && K_{16} \cdot c_4 c_5 .
\end{matrix} 
$$
Here $K_{13},K_{14},K_{15},K_{16} $ are cubic polynomials 
with $16$ terms indexed by trees  as in
(\ref{Kpolynomial}). By Proposition \ref{MTT},
they are the $3 \times 3$ minors of
the Laplacian of the complete graph ${\bf K}_4$:
$$ \begin{pmatrix}
\! \kappa_{13,14} \!+\! \kappa_{13,15} \!+\! \kappa_{13,16} \!\!
& -\kappa_{13,14} & -\kappa_{13,15} &- \kappa_{13,16} \\
 - \kappa_{14,13} 
& \!\! \kappa_{14,13}  \!+\!\kappa_{14,15} \!+\! \kappa_{14,16} \!\!
& - \kappa_{14,15} & - \kappa_{14,16} \\
-\kappa_{15,13} &- \kappa_{15,14} &
\!\! \kappa_{15,13} \!+\! \kappa_{15,14} \!+\! \kappa_{15,16} \!\!
& -\kappa_{15,16} \\
-\kappa_{16,13} & 
-\kappa_{16,14} & 
-\kappa_{16,15} & 
 \!\! \kappa_{16,13} \!+\! \kappa_{16,14} \!+\! \kappa_{16,15} 
\end{pmatrix}.
$$
The recombination network $G$ has $l = 7$ connected components
and its deficiency is $\delta = 5$, as there
$n=16$ complexes,
and the stoichiometric subspace $S$ has dimension  $\,\sigma = 4$.
The moduli ideal $M_G$  is minimally generated by $18$ binomials.
Twelve of them are cubics:
$$ \begin{matrix}
\,\, K_{8} K_{11} K_{15} - K_{7} K_{12} K_{16}\,\,\,\,&
\,\, K_{6} K_{9} K_{15} - K_{5} K_{10} K_{16}\, \, \,\, &
\,\, K_{4} K_{11} K_{14} - K_{3} K_{12} K_{16}\, \, \,\, \\ 
\,\, K_{2} K_{9} K_{14} - K_{1} K_{10} K_{16}\, \, \,\, &
\,\, K_{4} K_{7} K_{14} - K_{3} K_{8} K_{15}\, \, \,\, &
\,\, K_{2} K_{5} K_{14} - K_{1} K_{6} K_{15}\, \, \,\, \\
\,\, K_{6} K_{12} K_{13} - K_{5} K_{11} K_{14}\, \, \,\, &
\,\, K_{2} K_{12} K_{13} - K_{1} K_{11} K_{15}\, \, \,\, &
\,\, K_{8} K_{10} K_{13} - K_{7} K_{9} K_{14}\, \, \,\, \\
\,\, K_{4} K_{10} K_{13} - K_{3} K_{9} K_{15}\, \, \,\, &
\,\, K_{2} K_{8} K_{13} - K_{1} K_{7} K_{16}\, \, \,\, &
\,\, K_{4} K_{6} K_{13} - K_{3} K_{5} K_{16}.
\end{matrix}
$$
The remaining six generators of $M_G$ are quartics:
$$
\begin{matrix}
K_{9} K_{11} K_{14} K_{15} - K_{10} K_{12} K_{13} K_{16} \qquad &
K_{6} K_{8} K_{13} K_{15} - K_{5} K_{7} K_{14} K_{16} \\
K_{2} K_{4} K_{13} K_{14} - K_{1} K_{3} K_{15} K_{16} \qquad &
K_{5} K_{8} K_{10} K_{11} - K_{6} K_{7} K_{9} K_{12} \\
K_{1} K_{4} K_{10} K_{11} - K_{2} K_{3} K_{9} K_{12} \qquad &
K_{1} K_{4} K_{6} K_{7} - K_{2} K_{3} K_{5} K_{8}.
\end{matrix}
$$
The moduli space (of toric dynamical systems on $G$) is
the toric variety $V(M_G)$  defined by
these $18$ binomials. It has codimension $5$ and degree $56$.
For any recombination rates $\kappa^0 \in V_{>0}(M_G)$
and any starting point $c^0$ in the 
{\em population simplex} $ \Delta_7$,
the trajectory of 
the toric dynamical system (\ref{CRN}) 
stays in the $4$-dimensional polytope
$\,( c^0 + S) \,\cap \, \Delta_7\,$ and is conjectured to converge to
the Birch point $c^*$.
  \cite{Akin79} calls $c^*$ the {\em Wright point}. It generalizes the classical
{\em Hardy-Weinberg equilibrium} in the 2-locus system. \qed
\end{example}

\section{Detailed Balancing Systems} \label{sec:detailed}

In this section we discuss an important subclass of toric dynamical systems
called detailed balancing systems. Here, every edge of the digraph $G$
exists in both directions. We can thus identify $G = (V,E)$ with the 
underlying undirected graph $\tilde G  = (V,\tilde E)$,
where $\tilde E = \bigl\{ \{i,j\} \,: \,(i,j) \in E \bigr\}$.
For each undirected edge $\{i,j\} \in \tilde E$ of the graph $\tilde G$
we define an $n \times n$-matrix $A_\kappa^{\{i,j\}}$ as follows.
In rows $i,j$ and columns $i,j$ the matrix $A_\kappa^{\{i,j\}}$ equals
$$  \begin{pmatrix}
- \kappa_{ij} & \phantom{-}\kappa_{ij} \\
\phantom{-}\kappa_{ji} & -\kappa_{ji} 
\end{pmatrix},
$$
and all other entries of the matrix $A_\kappa^{\{i,j\}}$ are $0$. The
Laplacian of $G$ decomposes as 
\begin{equation}
\label{AKappaSum}
 A_\kappa \,\, = \, \sum_{\{i,j\} \in \tilde E} A_\kappa^{\{i,j\}} .
 \end{equation}
 A {\em detailed balancing system} is a dynamical system (\ref{CRN})
 for which the algebraic equations $\,\Psi(c) \cdot A_\kappa^{\{i,j\}} = 0\,$
for $\{i,j\} \in \tilde E$ admit a strictly positive solution $c^* \in \R^s_{>0}$.
In light of (\ref{AKappaSum}), every detailed balancing system is
a toric dynamical system, so the positive solution $c^*$
is unique and coincides with the Birch point.
As it is for toric dynamical systems,
the condition of being detailed balancing
 depends on the graph $\tilde G$ and the constants~$\kappa_{ij}$.

We rewrite this condition in terms of binomials
in  $\Q[c,\kappa]$. The two non-zero entries
of the row vector $\,\Psi(c) \cdot A_\kappa^{\{i,j\}}\,$
are $\,\kappa_{ij}c^{y_i} - \kappa_{ji}c^{y_j} \,$
and its negative. Moreover, we find
$$ 
 \Psi(c) \cdot A_\kappa^{\{i,j\}} \cdot Y \quad = \quad
 (\kappa_{ij}c^{y_i} - \kappa_{ji}c^{y_j})\cdot (y_j - y_i),
 $$
 and hence the right hand side of the dynamical system (\ref{CRN})
 can be rewritten as follows:
 \begin{equation} \label{AKappaSum2}
\Psi(c) \cdot A_\kappa \cdot Y  \,\,  = \,
 \sum_{\{i,j\} \in \tilde E} \!\!\  \Psi(c) \cdot A_\kappa^{\{i,j\}} \cdot Y
 \,\, = \,
 \sum_{\{i,j\} \in \tilde E} (\kappa_{ij}c^{y_i} - \kappa_{ji}c^{y_j}) \cdot (y_j - y_i).
  \end{equation}
  For a detailed balancing system,
  each summand in (\ref{AKappaSum2}) 
  vanishes at the Birch point $c^*$.
  
\begin{example} \label{Ex:K3again}
We revisit Example \ref{Ex:K3}.
Let $s=2, n=3$ and $\tilde G$ the complete graph on three
nodes labeled by $c_1^2, c_1 c_2 $ and $c_2^2$. The
dynamical system (\ref{specialCRN}) is now written as
$$ \frac{d}{dt} ( c_1,c_2) \, =\,
(\kappa_{12} c_1^2 - \kappa_{21} c_1 c_2) \cdot (-1,1) \,+\,
(\kappa_{13} c_1^2 - \kappa_{31} c_2^2) \cdot (-2,2) \,+\,
(\kappa_{23} c_1 c_2 - \kappa_{32} c_2^2) \cdot (-1,1) . $$
This is a detailed balancing system if and only if the following algebraic identities hold:
  \begin{equation}
  \label{moduli1again}
\kappa_{12}^2 \kappa_{31} - \kappa_{21}^2 \kappa_{13} \,\,\, = \,\,\,
\kappa_{23}^2 \kappa_{31} - \kappa_{32}^2 \kappa_{13} \,\, = \,\,\,
\kappa_{12} \kappa_{32} - \kappa_{21} \kappa_{23} \,\,\,= \,\,\, 0 .
\end{equation}
This defines a toric variety of codimension two which lies in
the hypersurface (\ref{moduli1}). \qed
\end{example}

\smallskip

To fit our discussion into the algebraic framework of Section 2,
we now propose the following definitions.
The {\em detailed balancing ideal} is the following toric ideal in $\Q[\kappa,c]$:
\begin{equation}
\label{DBideal}
 \widetilde T_G \,\,\, := \,\,\,
\bigl(\, \langle \, \kappa_{ij} c^{y_i} - \kappa_{ji} c^{y_j} \, \,| \,\,
\{i,j\} \in \tilde E \,\rangle \,\,: \,\, (c_1 c_2 \cdots c_s)^\infty \,\bigr). 
\end{equation}
The corresponding elimination ideal in $\Q[\kappa]$ will be called the
{\em detailed moduli ideal}:
$$ \widetilde M_G \quad := \quad \widetilde T_G \,\cap \,\Q[\kappa]. $$
The ideal $\widetilde T_G$ is toric, by the same reasoning as in 
Proposition  \ref{prop:itstoric}. The detailed moduli ideal
$\widetilde M_G$ is a toric ideal of Lawrence type, as was the ideal in
Example \ref{ex:lawrence}. Note, however, that the ideals
$\widetilde T_G$ and $\widetilde M_G$ are toric in the
original coordinates $\kappa_{ij}$. Here, we did not
need the transformation to the new
coordinates $K_1,\ldots,K_n$ in (\ref{Kpolynomial}).

Using the ring inclusion $\,\Q[K,c] \subset \Q[\kappa,c]$,
we have the following inclusions of ideals:
$$ T_G \,\subseteq \, \widetilde T_G \quad
\hbox{and} \quad M_G \,\subseteq \, \widetilde M_G .$$
Here the equality holds precisely in the situation of Example \ref{ex:lawrence},
namely, when each chemical complex appears in only one
reaction and each reaction is reversible.
In general, as seen in Example \ref{Ex:K3again},
the corresponding inclusion 
of moduli spaces will be strict:
$$ V_{>0} (\widetilde M_G) \,\,\subset \, \,V_{>0}(M_G). $$
In words: every detailed balancing system is a toric dynamical system
but not vice versa.

The following characterization of detailed balancing systems will
be used in the next section.
If $L$ is any vector in $\R^s$ 
and $c$ the unknown concentration vector  then we write
$$\, L*c \,\, := \,\, (L_1 c_1, L_2 c_2, \ldots, L_s c_s). $$

\begin{lem}
A toric dynamical system is detailed balancing if and only if all the binomials 
$\, \kappa_{ij}c^{y_i} - \kappa_{ji}c^{y_j}\,$
in (\ref{DBideal}) have the form $ (L * c)^{y_i} - (L*c)^{y_j}$, for some positive vector $L \in \R_{>0}^s$.
Thus, a detailed balancing system is a toric dynamical system 
of the special form
\begin{equation}
\label{DBS}
 \frac{d c}{dt} \quad = \quad \sum_{\{i,j\} \in \tilde E} \!
\bigl( (L * c)^{y_i} - (L*c)^{y_j} \bigr) \cdot (y_j - y_i).
 \end{equation}
\end{lem}

\begin{proof}
The if-direction is easy: if our binomials have the special form
$ (L * c)^{y_i} - (L*c)^{y_j}$ then
$\,c^* = (1/L_1,1/L_2,\ldots,1/L_s)\,$ is a positive
solution to the equations $\,\Psi(c) \cdot A_\kappa^{\{i,j\}} = 0$.
Conversely, for the only-if direction, we define $L$ as
the reciprocal of the Birch point
$\,L = (1/c_1^*, 1/c_2^*, \ldots, 1/c_n^*)$,
and the result follows the fact that $\,c^{y_i-y_j} = (c^*)^{y_i-y_j}\,$ remains valid for
all stationary points $c$ of the system (\ref{CRN}) as the starting point $c(0)$ varies.
\end{proof}

We now fix a detailed balancing system (\ref{DBS}) with a particular
starting point $c(0)$. Then the trajectory $c(t)$ evolves inside the
invariant polyhedron $\,  P \,=\,(c(0)+S) \,\cap \, \R^s_{\geq 0}$.
Consider any {\em acyclic orientation} $E' \subset \tilde E$ of the graph $\tilde G$. This means that
$E'$ contains  one from each pair of directed edges $(i,j)$ and $(j,i)$ in $E$,
in such a way that the resulting directed subgraph of $G$ has no directed cycles.
The acyclic orientation $E'$ specifies a {\em stratum} 
$\mathcal{S}$ inside the relatively open polyhedron
$\,  P^o \,=\,(c(0)+S) \,\cap \, \R^s_{> 0}\,$ as follows:
$$ \mathcal{S}
\,\,\,:= \,\,\, \bigl\{\, c \in P^o ~|~ (L*c )^{y_i}>(L*c)^{y_j} \mbox{ for  all $(i,j)$ in $E'$} \, \bigr\}.$$
The invariant polyhedron $P$ is partitioned into such strata
and their boundaries. We are interested
in how the strata meet the boundary of $P$. Each face 
 of $P$ has the form $\,F_I := \{ c \in P \, |\,
  c_i = 0 \mbox{ for } i \in I\}\, $ where $I $ is subset of $\{1,2,\ldots,s\}$.
This includes $F_\emptyset = P$.

 \begin{lem}\label{farkas}
 Consider a detailed balancing system (\ref{DBS}) and
 fix an acyclic orientation $E'$ of the
 graph $\tilde G$.   If the closure of the  stratum $\mathcal{S}$
 corresponding to $E'$
  intersects the relative interior of a face $ F_I$ of the invariant polyhedron $P$, 
  then there exists a strictly positive vector ${\alpha} \in \R^I_{>0}$ such that 
  $\,\,\sum_{k \in I } (y_{jk} - y_{ik}) \cdot {\alpha}_k \, \geq \,0\,\,$ for all 
directed edges  $(i,j)$ in $E'$.
\end{lem}

\begin{pf} 
We proceed by contradiction: assume that the inequalities 
$\,\sum_{k \in I } (y_{jk} - y_{ik}) {\alpha}_k \geq 0\,$ 
have no strictly positive solution $\alpha \in \R ^I _ {> 0}$.  
By Linear Programming Duality (Farkas' Lemma), there is a non-negative linear combination $v = \sum_{(i,j) \in E'} \lambda_{ij} (y_{j} - y_{i}) $ such that the following two conditions on $v$ hold: 
(a) supp$(v^+) \cap I = \emptyset$, and
(b) supp$(v^-)$ contains some $j_0 \in I$.  
We shall prove the following two claims, which give the
desired contradiction:

\smallskip

{\bf Claim One}: {\em If $c$ is a point in the relative interior of $F_I$, then $(L*c)^{v_+} >  (L*c) ^{v_-}$. }

\smallskip

Since $(L*c)_i =0$ if and only if $i \in I$, and $(L*c)_j > 0$ for all $j \notin I$,
(a) implies that $(L*c)^{v+}$ is strictly positive, while
(b) implies that $(L*c)^{v-} =0$, and we are done.

\smallskip

{\bf Claim Two}: {\em If $c $ is a point in the closure of the stratum $\mathcal{S}$, then
$(L*c)^{v_+}\leq (L*c) ^{v_-}$.}

\smallskip

Consider any point $s \in \mathcal{S}$. By the construction of $v$, the following equation holds:
\begin{equation}
\label{LSV}
(L*s)^v \,\, =  \,\, (L*s)^{  \sum_{(i,j) \in E'} \lambda_{ij} (y_{j} - y_{i}) }
	\,\, = \, \prod_{(i,j) \in E'} \left( (L*s)^{y_j - y_i} \right)^{\lambda _{ij}}.
\end{equation}
Recall that $\,(L*s)^{y_j - y_i} \leq 1\,$ for each oriented edge $(i,j)
\in E'$.   Also, each $\lambda_{ij}$ is non-negative, so  $
((L*s)^{y_j- y_i})^{\lambda _ {ij}} \leq 1$.  Using (\ref{LSV}), this implies that
$(L*s)^v \leq 1$, and therefore $(L*s)^{v_+} \leq (L*s)^{v_-}$.
By continuity we can replace $s$ by $c$ in this last inequality.
\end{pf}

The vector $\alpha \in \R_{>0}^I$ in Lemma~\ref{farkas} will play
a special role in the next section. In Corollary \ref{claim3} below
we regard $\alpha$ as a vector in $\R_{\geq 0}^s$
by setting $\alpha_j = 0$ for all $j \in \{1,\ldots,s\}\backslash I$.

\begin{cor}\label{claim3}
Let $c(t)$ be a trajectory of a detailed balancing system (\ref{DBS})
on the invariant polyhedron $P$,
and suppose that a point $c(t_0)$ on this trajectory lies both in the closure
of a stratum $\mathcal{S}$ and in the relative interior of a face $F_I$ of $P$.
 Let $\alpha \in \R^s_{\geq 0}$ be the vector obtained as in
 Lemma~\ref{farkas}. Then, the inner product $\,\langle \, \alpha, \,\frac{dc}{dt}(t_0) \,\rangle \,$
 is non-negative.
\end{cor}

\begin{pf}
Let $E' $ denote the orientation which specifies $\mathcal{S}$.  
The velocity vector
${\frac{dc}{dt}(t_0)}$ equals 
\[\sum_{(i,j) \in E'} \bigl( (L*c(t_0))^{y_i} - (L*c(t_0) )^{y_j}\bigr) \cdot (y_j - y_i). \]
Since $c(t_0)$ is in the closure of the stratum $\mathcal{S}$, we have
$\, (L * c(t_0))^{y_i} -  (L * c(t_0)) ^{y_j}  \, \geq \,0$.
We also have $\,\left\langle {\alpha}  ,  y_j-y_i \right\rangle \geq 0\,$ 
because ${\alpha}$ comes from Lemma \ref{farkas}.  This implies
\begin{align*}
\langle\, {\alpha}  \,, \ { \frac{dc}{dt}(t_0) }\, \rangle \,\,\,=\,\, \sum_{(i,j) \in E'}  \!
 \bigl( (L*c(t_0))^{y_i} - (L*c(t_0) )^{y_j}\bigr) \cdot \left\langle \, {\alpha} \, , \   y_j-y_i \ \right\rangle 
\, \,\,\geq \,\,\, 0.
\end{align*} 
This is the claimed inequality. It will be used in the proof of
Theorem \ref{MainConvergence}.
\end{pf}

\section{Partial Results on the Global Attractor Conjecture} \label{sec:GACresults}

This section contains what we presently know about the Global Attractor Conjecture 
which was stated in Section 3. This conjecture is proved
for detailed balancing systems
 whose invariant polyhedron is bounded and of dimension two. 
 We begin with  some general facts on trajectories
of toric dynamical systems, which are interesting in their own right. 

Consider a fixed toric dynamical system (\ref{CRN}) with
strictly positive starting point $c(0) = c^0  \in \R_{>0}^s$. The trajectory $c(t)$
remains in the invariant polyhedron $ P = (c^0 +S) \cap \R_{\geq 0}^s $.
Recall that any face of $P$ has the form $F_I:= \{ c \in P \, |\, 
  c_i = 0 \mbox{ if } i \in I\} $, where $I \subseteq \{1, \ldots s\}$.
      The boundary $\partial P$ of $P$ is the union of all faces $F_I$ 
    where $I$ is a proper subset of $\{1,\ldots,s\}$.
For positive $\varepsilon$, the $\varepsilon$-neighborhood in $P$ of the boundary of $P$ will be denoted by
$V_\varepsilon(\partial P)$. 

We note that the transformed entropy function~(\ref{eq:entropy})
can be extended continuously to the boundary of $P$, because $c_i \log c_i \to 0$ as $c_i\to 0^+$. 
Equivalent formulations of the following result are well known.
For instance, see \cite{SiegelChen94,Sontag01}.

\begin{prop} \label{prop:interior}
Suppose that the invariant polyhedron $P$ is bounded and 
the distance between the boundary of $P$ 
and the set $\{c(t)\in P \ | \ t>0\}$ is strictly positive. Then
the trajectory $c(t)$ converges to the Birch point $c^*$ of $P$.
\end{prop}

\begin{pf} 
We assume that $c(t)$ does not converge to $c^*$.
 Let $\varepsilon>0$ be such that $c(t) \notin V_\varepsilon(\partial P)$ for all $t>t_0$. 
The strict Lyapunov function~(\ref{eq:entropy}) ensures that there exists a 
neighborhood $V_{\varepsilon'}(c^*)$ of the Birch point $c^*$ such that all trajectories that visit 
$V_{\varepsilon'}(c^*)$ converge to $c^*$.
Then $c(t) \notin V_{\varepsilon'}(c^*)$ for all $t>t_0$.
Denote the complement of the two open neighborhoods by $P_0 := P \setminus \left(V_\varepsilon(\partial P) \cup V_{\varepsilon'}(c^*)\right)$. Then the non-positive and continuous function $c \mapsto (\nabla E \cdot \frac{dc}{dt})(c)$ does not vanish on $P_0$ by Proposition \ref{prop:Birch}, so it is bounded above by some $-\delta<0$ on $P_0$. Therefore, the value of $E(c(t))$ decreases at a rate of at least $\delta$ for all $t>t_0$, which implies that $E$ is unbounded on $P_0$.  This is a contradiction.
\end{pf}

Given a trajectory $c(t)$ of (\ref{CRN}),
a point $\bar c \in P$  is called an \emph{$\omega$-limit point}
 if  there exists a sequence $t_n \to \infty$ with $\lim_{n\to \infty} c(t_n) = \bar c$.
Proposition \ref{prop:interior}  says that if the trajectory $c(t)$ does not 
have any $\omega$-limit points on the boundary of $P$, then it must converge to the Birch point $c^*$.
Thus, in order to prove the Global Attractor Conjecture,
it would suffice to show that no boundary point of $P$ is an $\omega$-limit point.
We first rule out the vertices.

\begin{prop}\label{vertex}
Let $r$ be a vertex of $P$ and consider any $\varepsilon > 0$.
Then, there exists a neighborhood $W$ of $r$ such that  any trajectory 
$c(t)$ with starting point $c(0)=c^o$ satisfying  $dist(c^0,r)>\varepsilon$, 
does not visit $W$ for any $t>0$.
\end{prop}

\begin{pf}
The following set is the intersection of a closed cone with a sphere of radius one:
\begin{align*}
\mathcal{V} \,\, :=  \,\, \left\{ \frac{v}{\lVert v \rVert}  \ | \ v \in S \backslash \{0\}
\mbox{ and } r+ v \mbox{\, lies in \,} P \right\}.
\end{align*}
Hence $\mathcal{V}$ is compact.
We set $\,I \,=\, \bigl\{\, j \in \{1,\ldots,s\} \,: \, r_j = 0 \bigl\}$.
For each $v \in \mathcal{V}$, the ray 
$\,\gamma _ v (t) \,:=  r + t v \,$
extends from the vertex $\gamma_v (0) = r$ into the polyhedron $P$ for small $t>0$.  
We consider how the transformed entropy function changes along such a ray:
\begin{align*}
\frac{d}{dt} E(\gamma_v(t)) \,&=\,\, \sum_{j \in I } v_j  (\log (0 + t v_j) ) 
\,+\, \sum_{j \notin I } v_j \log (r_j + t v_j) \,-\, \sum_{i=1}^s \log ( c_j^* v_j)  \\
\,&= \,\, (\Sigma_{j \in I} v_j) \cdot \log (t)  \,\, + \,\, w(t),
\end{align*}
where the function $w(t)$ admits a universal upper bound for $t$ close to $0$ and $v \in \mathcal{V}$.
For each $j \in I$ we have $v_j \geq 0$ because $r_j = 0$ and $r + t v \in P$
for small $t > 0$. Also, since $v$ points into $P$, there exist $j \in I$
with $v_j > 0$. Thus, the function $\Sigma_{j \in I} v_j$ has a positive minimum over $\mathcal{V}$. It follows that $\frac{d} {dt} E (\gamma_v(t))$ tends to $- \infty$ for $t \to 0$. 
There exists $t_0 < \varepsilon$ such that for  all $v \in \mathcal{V}$
the function $\,t \mapsto E(r + t v)\,$ decreases for $ 0  < t \leq t_0$.
So, $E(r) > \mu:=\max_{v \in \mathcal{V}} E(r + t_0 v)$. On the other hand, $E$ is continuous, so there is a neighborhood $W$
of the vertex $r$ (contained in $\{ r + t v \,\,|\,\, t < t_0, v \in \mathcal{V} \}$) 
such that $\,E(c) >  (E(r)+ \mu)/2\,$ for all $c \in W$. Since $E$ decreases 
along trajectories, we conclude that no trajectory $c(t)$ that starts 
at distance $\geq \varepsilon$ from the vertex $r$ can enter $W$.
\end{pf}

\begin{remark} Chemical reaction networks for which $P$ is bounded are called \emph{conservative}.  For conservative
networks, there exists a positive mass assignment for each species that is conserved by all reactions \citep{Fein79}. On the other hand, if $0\in P$, then the reaction network is not conservative.  Proposition \ref{vertex} ensures that, for a toric dynamical system, complete depletion of all the concentrations $c_1, c_2, ..., c_s$ is never possible.
\end{remark}

\begin{lem}\label{claim1}
Suppose that $P$ is bounded and that the trajectory c(t) has an $\omega$-limit point on the boundary of $P$. Then for any $\varepsilon > 0$ there exists a positive number $t_\varepsilon >0$ such that $c(t)$ belongs to 
$V_\varepsilon(\partial P)$ for all $t > t_\varepsilon$.
In other words, the trajectory $c(t)$ approaches the boundary.
\end{lem}

\begin{pf}
Suppose that for some $\varepsilon > 0$ there exists a sequence $t_n \to \infty$ such that $c(t_n) \notin V_\varepsilon(\partial P)$ for all $n$. As $P$ is bounded, the trajectory $c(t)$ has an
$\omega$-limit point $p \in P \backslash V_\varepsilon(\partial P)$. 
On the other hand, $c(t)$ also has an $\omega$-limit point on the boundary of $P$. Consider a ball $B_{2\delta}(p)$ of radius $2\delta$ around $p$, whose closure lies fully in the relative interior of $P$. The trajectory $c(t)$ enters and exits the neighborhood $B_{\delta}(p)$ of $p$ infinitely many times, and also enters and exits the neighborhood $P \backslash B_{2\delta}(p)$ of the boundary infinitely many times. The trajectory $c(t)$ travels repeatedly between these two sets which are at distance $\delta$ from each other.   
Note that $|dc/dt|$ is bounded above, and $\nabla E\cdot dc/dt$ is bounded away from zero on the 
annulus $\,B_{2\delta}(p) \backslash B_{\delta}(p)$. Then, as in the proof of Proposition \ref{prop:interior}, each traversal between the neighborhoods decreases the value of $E(c(t))$ by a positive amount that is bounded away from zero.  This contradicts the fact that $E$ is bounded on $P$.
\end{pf}

 We shall now prove the main result of this section.
Admittedly, Theorem \ref{MainConvergence} has three
rather restrictive hypotheses, namely,
``dimension two,''  ``bounded polyhedron,'' and
``detailed balancing.''  At present we do not
know how to remove any of these hypotheses.

\begin{thm} \label{MainConvergence}
Consider a detailed balancing system (\ref{DBS}) whose
stoichiometric subspace $\,S = \R \{y_j-y_i \,|\, (i,j) \in \tilde E\}\,$
is two-dimensional and assume that
the invariant polygon $ \,P = (c^0 +S) \cap \R_{\geq 0}^s \,$ is bounded.
Then the Birch point $c^*$ is a global attractor for $P$. 
\end{thm}

\begin{pf}
By Proposition \ref{prop:interior}, we need only rule out the possibility that the trajectory $c(t)$ has an $\omega$-limit point on the boundary of $P$.  
Proposition \ref{vertex} gives the existence of open neighborhoods of the vertices such that no trajectory $c(t)$ that starts outside them can visit them. Let $V$ denote the union of these neighborhoods. 
Suppose now that $c(t)$ has an $\omega$-limit point on $\partial P$.
That limit point  lies in the relative interior of some edge $F$ of $P$.
Let $F_\varepsilon$ denote the set of points in $P$ which have distance
at most $\varepsilon$ from the edge $F$.

We  claim that there exists $\varepsilon > 0$ and $t_\varepsilon > 0$, 
such that the trajectory $c(t)$ remains in the subset $F_\varepsilon \backslash V$ 
for all $t > t_\varepsilon$. This is true because  $c(t)$ belongs to the neighborhood $V_\varepsilon (\partial P)$ of the boundary for $t \gg 0$,  by Lemma \ref{claim1},
and hence $c(t)$ belongs to $V_\varepsilon(\partial P)\backslash V$ for $t \gg 0$.
But this implies that $c(t)$ belongs to $F_\varepsilon \backslash V$ 
for $t \gg 0$ because $F_{\varepsilon} \backslash V$ 
is a connected component of $V_{\varepsilon}(\partial P) \backslash V$ for $\varepsilon$ sufficiently small. This uses the dimension two assumption.

Consider the closures of all strata $\mathcal{S}$ that intersect the relative interior of $F$.
After decreasing $\varepsilon$ if necessary, we may assume that
the union of these closures contains the set $F_\varepsilon \backslash V$,
which contains the trajectory $c(t)$ for $t > t_\varepsilon$. 
To complete the proof, we will show that
 the distance from $c(t)$ to the edge $F$ never decreases
 after $c(t)$ enters $\,F_\varepsilon \backslash V$.

Any stratum $\mathcal{S}$ whose closure intersects the relative interior of $F$
contributes a vector $\alpha = \alpha(\mathcal{S})$ which satisfies the statement of
Lemma~\ref{farkas} for $F = F_I$. The orthogonal projection of $\alpha(\mathcal{S})$ into
the two-dimensional stoichiometric subspace is a positive multiple
of the unit inner normal $\alpha_0 \in S$  to $F$ in $P$. 
By Corollary~\ref{claim3} we have 
$\, \langle \alpha(S), \frac{dc}{dt}(t) \rangle \geq 0 \,$ 
and hence $\,\langle \alpha_0,  \frac{dc}{dt}(t) \rangle \geq 0\,$
for $t  > t_\varepsilon$. Therefore the distance from $c(t)$ to $F$
cannot decrease. This is a contradiction to the assumption
that $F$ contains an $\omega$-limit point.
\end{pf}

\end{document}